# Projective varieties invariant by one-dimensional foliations

By Marcio G. Soares

*To Jacob Palis on his* 60<sup>th</sup> *anniversary*

## 1. Introduction

This work concerns the problem of relating characteristic numbers of one-dimensional holomorphic foliations of $\mathbb{P}^n_{\mathbb{C}}$ to those of algebraic varieties invariant by them. More precisely: if $M$ is a connected complex manifold, a one-dimensional holomorphic foliation $\mathcal{F}$ of $M$ is a morphism $\Phi : \mathcal{L} \longrightarrow \mathrm{T}M$ where $\mathcal{L}$ is a holomorphic line bundle on $M$. The singular set of $\mathcal{F}$ is the analytic subvariety $\mathrm{sing}(\mathcal{F}) = \{p : \Phi(p) = 0\}$ and the *leaves* of $\mathcal{F}$ are the leaves of the nonsingular foliation induced by $\mathcal{F}$ on $M \setminus \mathrm{sing}(\mathcal{F})$. If $M$ is $\mathbb{P}^n_{\mathbb{C}}$ then, since line bundles over $\mathbb{P}^n_{\mathbb{C}}$ are classified by the Chern class $c_1(\mathcal{L}) \in \mathrm{H}^2(\mathbb{P}^n_{\mathbb{C}}, \mathbb{Z}) \simeq \mathbb{Z}$, one-dimensional holomorphic foliations $\mathcal{F}$ of $\mathbb{P}^n_{\mathbb{C}}$ are given by morphisms $\Phi : \mathcal{O}(1-d) \longrightarrow \mathrm{T}\mathbb{P}^n_{\mathbb{C}}$ with $d \geq 0$, $d \in \mathbb{Z}$, which we call the *degree* of $\mathcal{F}$. We will use the notation $\mathcal{F}^d$ for such a foliation. Suppose now $V \xrightarrow{\boldsymbol{i}} \mathbb{P}^n_{\mathbb{C}}$ is an irreducible algebraic variety invariant by $\mathcal{F}^d$ in such a way that the pull-back $\boldsymbol{i}^* \mathcal{F}^d$ of $\mathcal{F}^d$ to $V$ has a finite set of points as the singular set. The problem we address is the relation between $d$ and the degree $d^0(V) = \int_V c_1(\boldsymbol{i}^* \mathcal{O}(1))^{\dim V}$ of $V$.

This question was first considered by Poincaré in [12], in the case of plane curves invariant by foliations of $\mathbb{P}^2_{\mathbb{C}}$ with a rational first integral. More recently Cerveau and Lins Neto [3], Carnicer [2] and Campillo and Carnicer [1] addressed this question when the invariant variety is a curve. In [14] we considered this problem for smooth hypersurfaces and in this work we treat the case of smooth algebraic varieties, concentrating on complete intersections, where effective computations are feasible.

To obtain the result we first calculate the number of singularities of the foliation in the invariant variety. To this end we impose the condition of nondegeneracy of the foliation along the variety; i.e., the linear part of a vector field representing the foliation on the variety has only nonzero eigenvalues at the singular points. With this at hand we can use Baum-Bott's theorem [4] and relate the degrees of the foliation and of the variety through a polynomial. Next we need a positivity argument saying that this polynomial is positive whenever it represents the number of singularities of a foliation leaving the variety invariant. We then try to obtain relations between the degrees. The positivity argument given here follows from the ampleness of the normal bundle of a smooth projective variety together with the vanishing theorem of [10] which, in turn, is a consequence of the work of Kamber and Tondeur [8] on foliated bundles.

## 2. Statement of results

Let $\mathcal{F}^d$ be a one-dimensional holomorphic foliation of $\mathbb{P}^n_\mathbb{C}$, $n \geq 2$, of degree $d \geq 2$, given by a morphism

$$\mathcal{O}(1-d) \xrightarrow{\Phi} T\mathbb{P}^n_\mathbb{C}$$

with singular set $\mathrm{sing}(\mathcal{F}^d) = \{p : \Phi(p) = 0\}$ which we assume to have codimension greater than 1. Such a foliation $\mathcal{F}^d$ is represented in affine coordinates $(z_1, \ldots, z_n)$ by a vector field of the form

$$X = gR + \sum_{j=0}^{d} X_j$$

where $R = \sum_{i=1}^{n} z_i \frac{\partial}{\partial z_i}$ is the radial vector field, $g \in \mathbb{C}[z_1, \ldots, z_n]$ is homogeneous of degree $d$ and $X_j$, $0 \leq j \leq d$, is a vector field whose components are homogeneous polynomials of degree $j$. Since $\mathrm{codim}\,\mathrm{sing}(\mathcal{F}^d) \geq 2$ we have either $g \not\equiv 0$ or $g \equiv 0$ and $X_d$ cannot be written as $hR$ with $h \in \mathbb{C}[z_1, \ldots, z_n]$ homogeneous of degree $d-1$. In this case $X$ has a pole of order $d-1$ along the hyperplane at infinity and it is worth remarking that, in case $g \not\equiv 0$, $g = 0$ is precisely the variety of tangencies of $\mathcal{F}^d$ with the hyperplane at infinity, whereas in case $g \equiv 0$ this hyperplane is invariant by the foliation.

If $i : V \longrightarrow \mathbb{P}^n_\mathbb{C}$ is a smooth algebraic variety invariant by $\mathcal{F}^d$, we say $\mathcal{F}^d$ is *nondegenerate along $V$* if its pull-back to $V$, $i^*\mathcal{F}^d$, is nondegenerate; i.e., $\mathrm{sing}(\mathcal{F}^d) \cap V$ is a finite set of points and, at all points $p \in \mathrm{sing}(\mathcal{F}^d) \cap V$, $\det\left(\frac{\partial Y_i}{\partial z_j}\right)(p) \neq 0$ for any vector field $Y$, tangent to $V$ at $p$ and representing $i^*\mathcal{F}^d$. Observe that we necessarily have $V \not\subset \mathrm{sing}(\mathcal{F}^d)$ and that a singularity of $\mathcal{F}^d$ which lies on $V$ may be degenerate, when considered as a singularity of $\mathcal{F}^d$ in $\mathbb{P}^n_\mathbb{C}$.

Let $\mathbf{V}^{n-k} \xrightarrow{i} \mathbb{P}^n_\mathbb{C}$, $n \geq 2$, be a smooth irreducible algebraic variety of codimension $k$. Define $\mathbf{V}^{n-k}_{[i]}$ as follows: $\mathbf{V}^{n-k}_{[0]} = \mathbf{V}^{n-k}$, $\mathbf{V}^{n-k}_{[1]}$ is a generic hyperplane section of $\mathbf{V}^{n-k}$ and $\mathbf{V}^{n-k}_{[i]}$ is a generic hyperplane section of $\mathbf{V}^{n-k}_{[i-1]}$, $i \geq 2$. In what follows it is irrelevant whether we regard $\mathbf{V}_{[i]}^{n-k}$ as a subvariety of $\mathbb{P}_{\mathbb{C}}^{n}$ or of $\mathbb{P}_{\mathbb{C}}^{n-i}$ and we let $\chi(\mathbf{V}_{[i]}^{n-k})$ denote the Euler-Poincaré characteristic of $\mathbf{V}_{[i]}^{n-k}$. Of course, by Bézout's theorem, $\chi(\mathbf{V}_{[n-k]}^{n-k})$ is just the degree $d^0(\mathbf{V}^{n-k})$ of $\mathbf{V}^{n-k}$.

We then have:

THEOREM I. *Let $\mathbf{V}^{n-k} \xrightarrow{i} \mathbb{P}_{\mathbb{C}}^{n}$ be a smooth irreducible algebraic variety invariant by $\mathcal{F}^d$, a one-dimensional holomorphic foliation of $\mathbb{P}_{\mathbb{C}}^{n}$ of degree $d \geq 2$, nondegenerate along $\mathbf{V}^{n-k}$. Then the number of singularities of $\mathcal{F}^d$ in $\mathbf{V}^{n-k}$, noted $\mathcal{N}(i^*\mathcal{F}^d, \mathbf{V}^{n-k})$, is*

$$\mathcal{N}(i^*\mathcal{F}^d, \mathbf{V}^{n-k}) = \chi(\mathbf{V}_{[n-k]}^{n-k})d^{n-k} + \sum_{j=1}^{n-k}\left[\chi(\mathbf{V}_{[n-k-j]}^{n-k}) - \chi(\mathbf{V}_{[n-k-j+1]}^{n-k})\right]d^{n-k-j}.$$

*Alternatively,*

$$\mathcal{N}(i^*\mathcal{F}^d, \mathbf{V}^{n-k}) = \sum_{j=0}^{n-k}\left[\sum_{i=0}^{j}(-1)^i \varrho_i(\mathbf{V}^{n-k})\right]d^{n-k-j}$$

*where $\varrho_i(\mathbf{V}^{n-k})$ is the $i^{\text{th}}$ class of $\mathbf{V}^{n-k}$. Moreover, $\mathcal{N}(i^*\mathcal{F}^d, \mathbf{V}^{n-k}) > 0$.*

Suppose now that $\mathbf{V}_{(d_1,\ldots,d_k)}^{n-k} \xrightarrow{i} \mathbb{P}_{\mathbb{C}}^{n}$ is a smooth irreducible complete intersection defined by $F_1 = 0, \ldots, F_k = 0$ where $F_\ell \in \mathbb{C}[z_0, \ldots, z_n]$ is homogeneous of degree $d_\ell$, $1 \leq \ell \leq k$. Then we have:

COROLLARY. *Assume $\mathbf{V}_{(d_1,\ldots,d_k)}^{n-k}$ is invariant by $\mathcal{F}^d$, a one-dimensional holomorphic foliation of $\mathbb{P}_{\mathbb{C}}^{n}$ of degree $d \geq 2$, nondegenerate along $\mathbf{V}_{(d_1,\ldots,d_k)}^{n-k}$. Then the number of singularities of $\mathcal{F}^d$ in $\mathbf{V}_{(d_1,\ldots,d_k)}^{n-k}$ is*

$$\mathcal{N}(i^*\mathcal{F}^d, \mathbf{V}_{(d_1,\ldots,d_k)}^{n-k}) = \sum_{j=0}^{n-k}\left[\sum_{\delta=0}^{j}(-1)^\delta \varrho_\delta\left(\mathbf{V}_{(d_1,\ldots,d_k)}^{n-k}\right)\right]d^{n-k-j}$$

$$= (d_1 \ldots d_k)\sum_{j=0}^{n-k}\left[\sum_{\delta=0}^{j}(-1)^\delta \mathcal{W}_\delta^{(k)}(d_1-1,\ldots,d_k-1)\right]d^{n-k-j}$$

*where $\mathcal{W}_\delta^{(k)}$ is the Wronski (or complete symmetric) function of degree $\delta$ in $k$ variables*

$$\mathcal{W}_\delta^{(k)}(X_1,\ldots,X_k) = \sum_{i_1+\cdots+i_k=\delta} X_1^{i_1}\ldots X_k^{i_k}.$$

As an application of Theorem I and the corollary we have:

THEOREM II. *Let* $\mathbf{V}^{n-k}_{(d_1,\ldots,d_k)}$ *and* $\mathcal{F}^d$ *be as in Theorem* I *and suppose* $n-k$ *odd. Then*

$$d \geq \frac{\varrho_{n-k}(\mathbf{V}^{n-k}_{(d_1,\ldots,d_k)})}{\varrho_{n-k-1}(\mathbf{V}^{n-k}_{(d_1,\ldots,d_k)})}.$$

In order to avoid trivialities, the invariant varieties considered in Theorem II are not linear subspaces of $\mathbb{P}^n_\mathbb{C}$. The theorem gives, in case $k = n-1$, that

$$d \geq \mathcal{W}_1^{(n-1)}(d_1 - 1, \ldots, d_{n-1} - 1),$$

and hence

$$d^0(\mathbf{V}^1_{(d_1,\ldots,d_{n-1})}) \leq \left(1 + \frac{d}{n-1}\right)^{n-1}.$$

Also, if $k = 1$ we obtain $d^0(\mathbf{V}^{n-1}_{d_1}) \leq d+1$, since $\varrho_j\left(\mathbf{V}^{n-1}_{d_1}\right) = d_1(d_1-1)^j$ (see [14]).

*Example* 1. Let $\mathbf{V}^{2n-1}_\ell$, $\ell \geq 3$, be the smooth hypersurface in $\mathbb{P}^{2n}_\mathbb{C}$ defined by $X_1^\ell + X_2^\ell + \cdots + X_{2n-1}^\ell + X_{2n}^\ell + X_{2n+1}^\ell = 0$. Also, $\mathbf{V}^{2n-1}_\ell$ is invariant by the foliation $\mathcal{F}$ of degree $\ell - 1$ on $\mathbb{P}^{2n}_\mathbb{C}$ defined by the vector field (affine coordinates $X_{2n+1} = 1$)

$$\mathcal{Z} = z_2^{\ell-1} z_1 \frac{\partial}{\partial z_1} + (z_2^\ell + 1) \frac{\partial}{\partial z_2}$$
$$+ \sum_{i=2}^n \left[ (z_2^{\ell-1} z_{2i-1} - z_{2i}^{\ell-1}) \frac{\partial}{\partial z_{2i-1}} + (z_2^{\ell-1} z_{2i} + z_{2i-1}^{\ell-1}) \frac{\partial}{\partial z_{2i}} \right].$$

Note that the bound $d^0(\mathbf{V}) = d+1$ is attained. Observe that in $\mathbb{P}^{2n}_\mathbb{C}$ a smooth hypersurface defines a hamiltonian vector field which can be used to foliate it.

*Example* 2. The elliptic quartic curve $\mathcal{E}_4$ can be realized as the complete intersection in $\mathbb{P}^3_\mathbb{C}$ defined by the quadrics $Q_1 = \{X_1^2 + X_2^2 + X_3^2 + X_4^2 = 0\}$ and $Q_2 = \{X_1 X_3 + X_2 X_4 = 0\}$. The foliation $\mathcal{F}$ of degree 2, on $\mathbb{P}^3_\mathbb{C}$, defined by the vector field (affine coordinates $X_4 = 1$)

$$\mathcal{Z} = (-z_1^2 z_2 + z_1 z_3) \frac{\partial}{\partial z_1} + (-z_1 z_2^2 + 2 z_2 z_3 - z_1) \frac{\partial}{\partial z_2} + (-z_1 z_2 z_3 - z_2^2 + z_3^2 + 1) \frac{\partial}{\partial z_3}$$

has $\mathcal{E}_4$ as invariant curve. Note that the bound $d^0(\mathbf{V}) = \left(1 + \frac{d}{n-1}\right)^{n-1}$ is attained.

## 3. Proof of Theorem I

Choose a hyperplane $H_\infty$ transverse to $\mathbf{V}^{n-k}$ and such that $H_\infty \cap \mathbf{V}^{n-k} \cap \text{sing}(i^*\mathcal{F}^d) = \emptyset$. Since a vector field $X$ representing $\mathcal{F}^d$ has a pole of order $d-1$ along $H_\infty$, the same holds for the representative $X_{|\mathbf{V}^{n-k}}$ so that it defines a section of $\mathrm{T}\mathbf{V}^{n-k} \otimes i^*\mathcal{O}(d-1)$. This section has isolated nondegenerate zeros by hypothesis, and so, according to Baum-Bott's theorem [4] applied to the top Chern class:

$$\mathcal{N}(i^*\mathcal{F}^d, \mathbf{V}^{n-k}) = \int_{\mathbf{V}^{n-k}} c_{n-k}(\mathrm{T}\mathbf{V}^{n-k} \otimes i^*\mathcal{O}(d-1))$$

where integration is over the fundamental class of $\mathbf{V}^{n-k}$. Since

$$c_{n-k}(\mathrm{T}\mathbf{V}^{n-k} \otimes i^*\mathcal{O}(d-1))$$
$$= c_{n-k}(\mathrm{T}\mathbf{V}^{n-k} \otimes i^*\mathcal{O}(-1) \otimes i^*\mathcal{O}(d))$$
$$= \sum_{j=0}^{n-k} c_j(\mathrm{T}\mathbf{V}^{n-k} \otimes i^*\mathcal{O}(-1)) c_1(i^*\mathcal{O}(d))^{n-k-j}$$
$$= \sum_{j=0}^{n-k} c_j(\mathrm{T}\mathbf{V}^{n-k} \otimes i^*\mathcal{O}(-1)) c_1(i^*\mathcal{O}(1))^{n-k-j} d^{n-k-j},$$

and

$$c_j(\mathrm{T}\mathbf{V}^{n-k} \otimes i^*\mathcal{O}(-1)) = \sum_{i=0}^{j} (-1)^{j-i} \binom{n-k-i}{j-i} c_i(\mathbf{V}^{n-k}) c_1(i^*\mathcal{O}(1))^{j-i},$$

we get

(1)
$$c_{n-k}(\mathrm{T}\mathbf{V}^{n-k} \otimes i^*\mathcal{O}(d-1))$$
$$= \sum_{j=0}^{n-k} \left[ \sum_{i=0}^{j} (-1)^{j-i} \binom{n-k-i}{j-i} c_i(\mathbf{V}^{n-k}) c_1(i^*\mathcal{O}(1))^{j-i} \right] c_1(i^*\mathcal{O}(1))^{n-k-j} d^{n-k-j}$$
$$= \sum_{j=0}^{n-k} \left[ \sum_{i=0}^{j} (-1)^{j-i} \binom{n-k-i}{j-i} c_i(\mathbf{V}^{n-k}) c_1(i^*\mathcal{O}(1))^{n-k-i} \right] d^{n-k-j}.$$

Following Fulton [5, 14.4.15], recall that the cycle associated to the $j^{\text{th}}$ polar locus of $\mathbf{V}^{n-k}$ is given, for a general linear subspace $L^{k+j-2}$ of $\mathbb{P}_{\mathbb{C}}^n$, by

$$[\mathbf{V}^{n-k}(L^{k+j-2})] = \sum_{i=0}^{j} (-1)^i \binom{n-k+1-i}{j-i} c_i(\mathbf{V}^{n-k}) c_1(i^*\mathcal{O}(1))^{j-i}$$

and that the $j^{\text{th}}$ class $\varrho_j$ of $\mathbf{V}^{n-k}$ is defined to be the degree of $[\mathbf{V}^{n-k}(L^{k+j-2})]$, so that

$$\varrho_j = \int_{\mathbf{V}^{n-k}} \sum_{i=0}^{j} (-1)^i \binom{n-k+1-i}{j-i} c_i(\mathbf{V}^{n-k}) c_1(i^*\mathcal{O}(1))^{n-k-i}$$

since the degree is computed through multiplication by $c_1(i^*\mathcal{O}(1))^{n-k-j}$. Now, using Stifel's relation $\binom{n-k-i+1}{\ell-i} = \binom{n-k-i}{\ell-i} + \binom{n-k-i}{\ell-i-1}$ we get

$$\sum_{i=0}^{j} (-1)^i \varrho_i = \int_{\mathbf{V}^{n-k}} \sum_{i=0}^{j} (-1)^{j-i} \binom{n-k-i}{j-i} c_i(\mathbf{V}^{n-k}) c_1(i^*\mathcal{O}(1))^{n-k-i}.$$

It follows from (1) that

$$\mathcal{N}(i^*\mathcal{F}^d, \mathbf{V}^{n-k}) = \int_{\mathbf{V}^{n-k}} c_{n-k}(\mathrm{T}\mathbf{V}^{n-k} \otimes i^*\mathcal{O}(d-1))$$

$$= \sum_{j=0}^{n-k} \left[\sum_{i=0}^{j} (-1)^i \varrho_i\right] d^{n-k-j}.$$

Let us now recall a consequence of Lefschetz' theorem on hyperplane sections [9]. If $X$ is a smooth irreducible algebraic variety and $\mathcal{H}_{t \in \mathbb{P}^1_{\mathbb{C}}}$ is a pencil of hyperplanes with axis $L^{n-2}$ then the Euler-Poincaré characteristics are related by

$$\chi(X) = 2\chi(X \cap H) - \chi(X \cap L^{n-2}) + (-1)^{\dim X} \varrho_{\dim X}(X)$$

where $H$ is a generic element of the pencil. Applying this to $\mathbf{V}^{n-k}$ we get:

$$\chi(\mathbf{V}^{n-k}) - \chi(\mathbf{V}^{n-k}_{[1]}) = \chi(\mathbf{V}^{n-k}_{[1]}) - \chi(\mathbf{V}^{n-k}_{[2]}) + (-1)^{n-k} \varrho_{n-k}(\mathbf{V}^{n-k}).$$

By repeating this argument, using the Piene-Severi comparison theorem [11] (which says that the class of a hyperplane section of $X$ is $\varrho_{\dim X - 1}(X)$) we obtain

$$\sum_{i=0}^{j} (-1)^i \varrho_i(\mathbf{V}^{n-k}) = \chi(\mathbf{V}^{n-k}_{[n-k-j]}) - \chi(\mathbf{V}^{n-k}_{[n-k-j+1]})$$

so that

$$\mathcal{N}(i^*\mathcal{F}^d, \mathbf{V}^{n-k}) = \chi(\mathbf{V}^{n-k}_{[n-k]}) d^{n-k} + \sum_{j=1}^{n-k} \left[\chi(\mathbf{V}^{n-k}_{[n-k-j]}) - \chi(\mathbf{V}^{n-k}_{[n-k-j+1]})\right] d^{n-k-j}.$$

It remains to show $\mathcal{N}(i^*\mathcal{F}^d, \mathbf{V}^{n-k}) > 0$. To this end we invoke the vanishing theorem of [10, théorème 2] which states that, if $\mathbf{V}^{n-k}$ is foliated by $i^*\mathcal{F}^d$ without singularities, then any polynomial on the Chern classes of the normal bundle $\mathrm{N}_{\mathbf{V}^{n-k}}$ of $\mathbf{V}^{n-k}$ in $\mathbb{P}^n_{\mathbb{C}}$ must vanish in dimension greater than $2s$, where $s$ is the complex codimension of $i^*\mathcal{F}^d$. Now, codim $i^*\mathcal{F}^d = n-k-1$ and, since det $(N_{\mathbf{V}^{n-k}})$ is ample [6], it follows from the hard Lefschetz theorem [9] that the rational class $c_1(\det(N_{\mathbf{V}^{n-k}}))^{n-k}$ is a basis of $H^{2n-2k}(\mathbf{V}^{n-k},\mathbb{Q})$, and therefore $c_1(N_{\mathbf{V}^{n-k}})^{n-k}$ is nonzero. This finishes the proof of Theorem I.

## 4. Proof of the corollary to Theorem I

Let us calculate $c_j(\mathbf{V}^{n-k}_{(d_1,\ldots,d_k)})$. Set $c_1(i^*\mathcal{O}(1)) = \mathbf{h}$. It is well known [7] that the total Chern class of $\mathbf{V}^{n-k}_{(d_1,\ldots,d_k)}$ is given by

$$(2) \quad c(\mathbf{V}^{n-k}_{(d_1,\ldots,d_k)}) = \frac{(1+\mathbf{h})^{n+1}}{\prod_{\ell=1}^{k}(1+d_\ell\mathbf{h})}.$$

Recalling that the Wronski functions are defined by

$$\frac{1}{\prod_{\ell=1}^{k}(1+d_\ell\mathbf{t})} = \sum_{\delta=0}^{\infty}(-1)^\delta \, \mathcal{W}^{(k)}_\delta(d_1,\ldots,d_k) \, \mathbf{t}^\delta$$

we have that (2) becomes

$$c(\mathbf{V}^{n-k}_{(d_1,\ldots,d_k)}) = \sum_{j=0}^{n-k}\left[\sum_{i+\delta=j}\binom{n+1}{i}(-1)^\delta \, \mathcal{W}^{(k)}_\delta\right]\mathbf{h}^j$$

$$= \sum_{j=0}^{n-k}\left[\sum_{\delta=0}^{j}(-1)^\delta\binom{n+1}{j-\delta}\mathcal{W}^{(k)}_\delta\right]\mathbf{h}^j$$

so that

$$c_i(\mathbf{V}^{n-k}_{(d_1,\ldots,d_k)}) = \left[\sum_{\delta=0}^{i}(-1)^\delta\binom{n+1}{i-\delta}\mathcal{W}^{(k)}_\delta(d_1,\ldots,d_k)\right]\mathbf{h}^i, \quad 0 \le i \le n-k.$$

Hence,

$$(3) \quad \chi(\mathbf{V}^{n-k}_{(d_1,\ldots,d_k)}) = \int_{\mathbf{V}^{n-k}_{(d_1,\ldots,d_k)}} c_{n-k}(\mathbf{V}^{n-k}_{(d_1,\ldots,d_k)})$$

$$= \left[\sum_{\delta=0}^{n-k}(-1)^\delta\binom{n+1}{n-k-\delta}\mathcal{W}^{(k)}_\delta(d_1,\ldots,d_k)\right]\int_{\mathbf{V}^{n-k}_{(d_1,\ldots,d_k)}}\mathbf{h}^{n-k}$$

$$= (d_1\cdots d_k)\left[\sum_{\delta=0}^{n-k}(-1)^\delta\binom{n+1}{n-k-\delta}\mathcal{W}^{(k)}_\delta(d_1,\ldots,d_k)\right].$$

To calculate the Euler-Poincaré characteristic of the variety $\mathbf{V}^{n-k-(q)}_{(d_1,\ldots,d_k,1^q)}$, obtained by cutting $\mathbf{V}^{n-k}_{(d_1,\ldots,d_k)}$ successively by $q$ generic hyperplanes, we either add $q$ extra equations of degree 1, or regard $\mathbf{V}^{n-k-(q)}_{(d_1,\ldots,d_k,1^q)}$ as a complete intersection in $\mathbf{P}^{n-q}_{\mathbf{C}}$, given by $k$ equations of degrees $d_1,\ldots,d_k$. Doing it this last way we have, from (3):

$$\chi\bigl(\mathbf{V}^{n-k-(q)}_{(d_1,\ldots,d_k,1^q)}\bigr) = (d_1\cdots d_k)\left[\sum_{\delta=0}^{n-k-q}(-1)^\delta\binom{n-q+1}{n-q-k-\delta}\mathcal{W}^{(k)}_\delta(d_1,\ldots,d_k)\right]$$

and

$$\chi\bigl(\mathbf{V}^{n-k-(q+1)}_{(d_1,\ldots,d_k,1^{q+1})}\bigr)$$
$$= (d_1\cdots d_k)\left[\sum_{\delta=0}^{n-k-q-1}(-1)^\delta\binom{n-q}{n-q-k-1-\delta}\mathcal{W}^{(k)}_\delta(d_1,\ldots,d_k)\right].$$

The coefficient of $d^q$ in the formula of Theorem I is (by Stifel's relation):

$$\chi\bigl(\mathbf{V}^{n-k-(q)}_{(d_1,\ldots,d_k,1^q)}\bigr) - \chi\bigl(\mathbf{V}^{n-k-(q+1)}_{(d_1,\ldots,d_k,1^{q+1})}\bigr)$$
$$= (d_1\cdots d_k)\left[\sum_{\delta=0}^{n-k-q}(-1)^\delta\binom{n-q}{n-q-k-\delta}\mathcal{W}^{(k)}_\delta(d_1,\ldots,d_k)\right].$$

Setting $q = n-k-j$ we have that Theorem I reads

$$\mathcal{N}\bigl(i^*\mathcal{F}^d,\mathbf{V}^{n-k}_{(d_1,\ldots,d_k)}\bigr)$$
$$= (d_1\cdots d_k)\sum_{j=0}^{n-k}\left[\sum_{\delta=0}^{j}(-1)^\delta\binom{k+j}{j-\delta}\mathcal{W}^{(k)}_\delta(d_1,\ldots,d_k)\right]d^{n-k-j}.$$

Now, Lemma 2 of Todd ([15, p. 200]) tells us that

$$\mathcal{W}^{(k)}_p(d_1-1,\ldots,d_k-1) = \sum_{i=0}^{p}(-1)^{p-i}\binom{k+p-1}{p-i}\mathcal{W}^{(k)}_i(d_1,\ldots,d_k).$$

Taking the alternate sum and using Stifel's relation, we arrive at:

$$\sum_{\delta=0}^{j}(-1)^\delta\,\mathcal{W}^{(k)}_\delta(d_1-1,\ldots,d_k-1) = \sum_{\delta=0}^{j}(-1)^\delta\binom{k+j}{j-\delta}\mathcal{W}^{(k)}_\delta(d_1,\ldots,d_k).$$

Thus we recover the classical formulas of Severi [13] and Todd [15], for the classes of a smooth complete intersection:

$$\varrho_j\bigl(\mathbf{V}^{n-k}_{(d_1,\ldots,d_k)}\bigr) = (d_1\cdots d_k)\mathcal{W}^{(k)}_j(d_1-1,\ldots,d_k-1)\ ,\quad 0\le j\le n-k.$$

This finishes the proof of the corollary.

## 5. Proof of Theorem II

LEMMA 1. *Let $x_1, \ldots, x_k$ be nonnegative integers. Then*

$$\mathcal{W}_1^{(k)}(x_1, \ldots, x_k) \, \mathcal{W}_{\delta-1}^{(k)}(x_1, \ldots, x_k) - \mathcal{W}_\delta^{(k)}(x_1, \ldots, x_k)$$
$$\geq \mathcal{W}_1^{(k)}(x_1, \ldots, x_k) \, \mathcal{W}_{\delta-2}^{(k)}(x_1, \ldots, x_k) - \mathcal{W}_{\delta-1}^{(k)}(x_1, \ldots, x_k).$$

*Proof.* Just observe that every monomial appearing in $\mathcal{W}_i^{(k)}(x_1, \ldots, x_k)$ also appears in $\mathcal{W}_1^{(k)}(x_1, \ldots, x_k) \, \mathcal{W}_{i-1}^{(k)}(x_1, \ldots, x_k)$ with coefficient at least 1. In particular both sides of the inequality are nonnegative. Now, the left side is a sum of monomials of degree $\delta$ and the right side is a sum of a smaller or equal number of monomials of degree $\delta - 1$. Since $x_1, \ldots, x_k$ are nonnegative integers the result follows.

Note that Lemma 1 gives

$$\mathcal{W}_1^{(k)}(d_1 - 1, \ldots, d_k - 1)$$
$$\geq \frac{\mathcal{W}_\delta^{(k)}(d_1 - 1, \ldots, d_k - 1) - \mathcal{W}_{\delta-1}^{(k)}(d_1 - 1, \ldots, d_k - 1)}{\mathcal{W}_{\delta-1}^{(k)}(d_1 - 1, \ldots, d_k - 1) - \mathcal{W}_{\delta-2}^{(k)}(d_1 - 1, \ldots, d_k - 1)}.$$

LEMMA 2. *Suppose $n - k$ odd and $1 \leq k < n - 1$. If*

$$d < \min_{2 \leq \delta \leq n-k} \left\{ \frac{\mathcal{W}_\delta^{(k)}(d_1 - 1, \ldots, d_k - 1) - \mathcal{W}_{\delta-1}^{(k)}(d_1 - 1, \ldots, d_k - 1)}{\mathcal{W}_{\delta-1}^{(k)}(d_1 - 1, \ldots, d_k - 1) - \mathcal{W}_{\delta-2}^{(k)}(d_1 - 1, \ldots, d_k - 1)} \right\}$$

*then $\mathcal{N}(i^*\mathcal{F}^d, \mathbf{V}_{(d_1, \ldots, d_k)}^{n-k}) \leq 0$.*

*Proof.* To avoid cumbersome notation write $\mathcal{W}_\delta^{(k)}(d_1 - 1, \ldots, d_k - 1)$ as $\mathcal{W}_\delta^{(k)}$. By the corollary

$$(d_1 \cdots d_k)^{-1} \mathcal{N}(i^*\mathcal{F}^d, \mathbf{V}_{(d_1, \ldots, d_k)}^{n-k}) = d^{n-k} + \left[1 - \mathcal{W}_1^{(k)}\right] d^{n-k-1}$$
$$+ \left[1 - \mathcal{W}_1^{(k)} + \mathcal{W}_2^{(k)}\right] d^{n-k-2} + \left[1 - \mathcal{W}_1^{(k)} + \mathcal{W}_2^{(k)} - \mathcal{W}_3^{(k)}\right] d^{n-k-3} + \cdots$$
$$+ \left[1 - \mathcal{W}_1^{(k)} + \cdots + (-1)^j \mathcal{W}_j^{(k)}\right] d^{n-k-j}$$
$$+ \left[1 - \mathcal{W}_1^{(k)} + \cdots + (-1)^{j+1} \mathcal{W}_{j+1}^{(k)}\right] d^{n-k-j-1}$$
$$+ \cdots + \left[1 - \mathcal{W}_1^{(k)} + \cdots + \mathcal{W}_{n-k-1}^{(k)}\right] d + \left[1 - \mathcal{W}_1^{(k)} + \cdots - \mathcal{W}_{n-k}^{(k)}\right].$$

Grouping the terms pairwise, always assuming the term of highest degree in $d$ to be odd we get:

$$(d_1 \cdots d_k)^{-1} \mathcal{N}(i^*\mathcal{F}^d, \mathbf{V}^{n-k}_{(d_1,\ldots,d_k)}) = \left[d + 1 - \mathcal{W}_1^{(k)}\right] d^{n-k-1}$$
$$+ \left[d(1 - \mathcal{W}_1^{(k)} + \mathcal{W}_2^{(k)}) + (1 - \mathcal{W}_1^{(k)} + \mathcal{W}_2^{(k)} - \mathcal{W}_3^{(k)})\right] d^{n-k-3} + \cdots$$
$$+ \left[d(1 - \mathcal{W}_1^{(k)} + \cdots + \mathcal{W}_j^{(k)}) + (1 - \mathcal{W}_1^{(k)} + \cdots + \mathcal{W}_j^{(k)} - \mathcal{W}_{j+1}^{(k)})\right] d^{n-k-j-1}$$
$$+ \cdots + \left[d(1 - \mathcal{W}_1^{(k)} + \cdots + \mathcal{W}_{n-k-1}^{(k)}) + (1 - \mathcal{W}_1^{(k)} + \cdots + \mathcal{W}_{n-k-1}^{(k)} - \mathcal{W}_{n-k}^{(k)})\right].$$

The term preceeding $d^{n-k-j-1}$ can be regrouped as

$$(d + 1 - \mathcal{W}_1^{(k)}) + (-\mathcal{W}_1^{(k)}d + \mathcal{W}_2^{(k)} + \mathcal{W}_2^{(k)}d - \mathcal{W}_3^{(k)})$$
$$+ \cdots + (-\mathcal{W}_{j-1}^{(k)}d + \mathcal{W}_j^{(k)} + \mathcal{W}_j^{(k)}d - \mathcal{W}_{j+1}^{(k)}).$$

Now, Lemma 1 and the hypothesis imply

$$d + 1 - \mathcal{W}_1^{(k)} \leq 0$$

and

$$-\mathcal{W}_{j-1}^{(k)}d + \mathcal{W}_j^{(k)} + \mathcal{W}_j^{(k)}d - \mathcal{W}_{j+1}^{(k)} < 0,$$

so that $\mathcal{N}(i^*\mathcal{F}^d, \mathbf{V}^{n-k}_{(d_1,\ldots,d_k)}) \leq 0$. This proves the lemma.

LEMMA 3. *Let $x_1, \ldots, x_k$ be nonnegative integers. Then, for $j \geq 1$*

$$\left(\mathcal{W}_j^{(k)}\right)^2(x_1, \ldots, x_k) \geq \mathcal{W}_{j-1}^{(k)}(x_1, \ldots, x_k) \mathcal{W}_{j+1}^{(k)}(x_1, \ldots, x_k).$$

*Also,*

$$\min_{1 \leq j \leq n-k} \left\{ \frac{\mathcal{W}_j^{(k)}(x_1, \ldots, x_k)}{\mathcal{W}_{j-1}^{(k)}(x_1, \ldots, x_k)} \right\} = \frac{\mathcal{W}_{n-k}^{(k)}(x_1, \ldots, x_k)}{\mathcal{W}_{n-k-1}^{(k)}(x_1, \ldots, x_k)}.$$

*Proof.* Let us write $\mathcal{W}_j^{(m)}(x_1, \ldots, x_m)$ as $\mathcal{W}_j^{(m)}$. The proof is by induction on the number of variables. If $k = 1$ then

$$\left(\mathcal{W}_j^{(1)}\right)^2(x_1) = x_1^{2j} \geq x_1^{j-1} x_1^{j+1} = \mathcal{W}_{j-1}^{(1)}(x_1) \mathcal{W}_{j+1}^{(1)}(x_1).$$

Assume it holds for $k-1$, so that $\left(\mathcal{W}_j^{(k-1)}\right)^2 \geq \mathcal{W}_{j-1}^{(k-1)} \mathcal{W}_{j+1}^{(k-1)}$. Observe that this inequality implies

$$(*) \qquad \frac{\mathcal{W}_1^{(k-1)}}{\mathcal{W}_0^{(k-1)}} \geq \frac{\mathcal{W}_2^{(k-1)}}{\mathcal{W}_1^{(k-1)}} \geq \cdots \geq \frac{\mathcal{W}_j^{(k-1)}}{\mathcal{W}_{j-1}^{(k-1)}} \geq \frac{\mathcal{W}_{j+1}^{(k-1)}}{\mathcal{W}_j^{(k-1)}} \geq \cdots \geq \frac{\mathcal{W}_{n-k}^{(k-1)}}{\mathcal{W}_{n-k-1}^{(k-1)}}.$$

Now note that, since

$$(**) \quad \mathcal{W}_{j-1}^{(k)} = \mathcal{W}_0^{(k-1)} x_k^{j-1} + \mathcal{W}_1^{(k-1)} x_k^{j-2} + \cdots + \mathcal{W}_{j-1}^{(k-1)}$$

we have

$$\mathcal{W}_j^{(k)} = \mathcal{W}_0^{(k-1)} x_k^j + \mathcal{W}_1^{(k-1)} x_k^{j-1} + \cdots + \mathcal{W}_{j-1}^{(k-1)} x_k + \mathcal{W}_j^{(k-1)}$$
$$= \mathcal{W}_{j-1}^{(k)} x_k + \mathcal{W}_j^{(k-1)}$$

and

$$\mathcal{W}_{j+1}^{(k)} = \mathcal{W}_0^{(k-1)} x_k^{j+1} + \mathcal{W}_1^{(k-1)} x_k^j + \cdots + \mathcal{W}_{j-1}^{(k-1)} x_k^2 + \mathcal{W}_j^{(k-1)} x_k + \mathcal{W}_{j+1}^{(k-1)}$$
$$= \mathcal{W}_{j-1}^{(k)} x_k^2 + \mathcal{W}_j^{(k-1)} x_k + \mathcal{W}_{j+1}^{(k-1)}.$$

With this at hand we get

$$\left(\mathcal{W}_j^{(k)}\right)^2 = \left(\mathcal{W}_{j-1}^{(k)}\right)^2 x_k^2 + 2\mathcal{W}_{j-1}^{(k)} \mathcal{W}_j^{(k-1)} x_k + \left(\mathcal{W}_j^{(k-1)}\right)^2$$

and

$$\mathcal{W}_{j-1}^{(k)} \mathcal{W}_{j+1}^{(k)} = \left(\mathcal{W}_{j-1}^{(k)}\right)^2 x_k^2 + \mathcal{W}_{j-1}^{(k)} \mathcal{W}_j^{(k-1)} x_k + \mathcal{W}_{j-1}^{(k)} \mathcal{W}_{j+1}^{(k-1)}.$$

Hence,

$$(***) \quad \left(\mathcal{W}_j^{(k)}\right)^2 - \mathcal{W}_{j-1}^{(k)} \mathcal{W}_{j+1}^{(k)} = \left(\mathcal{W}_j^{(k-1)}\right)^2 + \mathcal{W}_{j-1}^{(k)} \mathcal{W}_j^{(k-1)} x_k - \mathcal{W}_{j-1}^{(k)} \mathcal{W}_{j+1}^{(k-1)}.$$

Let us consider the term

$$\mathcal{W}_{j-1}^{(k)} \mathcal{W}_j^{(k-1)} x_k - \mathcal{W}_{j-1}^{(k)} \mathcal{W}_{j+1}^{(k-1)} = \mathcal{W}_{j-1}^{(k)} \left(\mathcal{W}_j^{(k-1)} x_k - \mathcal{W}_{j+1}^{(k-1)}\right).$$

Using $(**)$ we obtain

$$\mathcal{W}_{j-1}^{(k)} \left(\mathcal{W}_j^{(k-1)} x_k - \mathcal{W}_{j+1}^{(k-1)}\right)$$
$$= \mathcal{W}_j^{(k-1)} x_k^j + \left(\mathcal{W}_1^{(k-1)} \mathcal{W}_j^{(k-1)} - \mathcal{W}_0^{(k-1)} \mathcal{W}_{j+1}^{(k-1)}\right) x_k^{j-1}$$
$$+ \left(\mathcal{W}_2^{(k-1)} \mathcal{W}_j^{(k-1)} - \mathcal{W}_1^{(k-1)} \mathcal{W}_{j+1}^{(k-1)}\right) x_k^{j-2} + \cdots$$
$$+ \left(\mathcal{W}_{j-1}^{(k-1)} \mathcal{W}_j^{(k-1)} - \mathcal{W}_{j-2}^{(k-1)} \mathcal{W}_{j+1}^{(k-1)}\right) x_k - \mathcal{W}_{j-1}^{(k-1)} \mathcal{W}_{j+1}^{(k-1)}.$$

By $(*)$, all the coefficients of $x_k^\ell$ are nonnegative, for $\ell \geq 1$. Taking this into $(***)$ we conclude that

$$\left(\mathcal{W}_j^{(k)}\right)^2 - \mathcal{W}_{j-1}^{(k)} \mathcal{W}_{j+1}^{(k)} \geq \left(\mathcal{W}_j^{(k-1)}\right)^2 - \mathcal{W}_{j-1}^{(k-1)} \mathcal{W}_{j+1}^{(k-1)} \geq 0$$

by inductive hypothesis. From this it follows that

$$\min_{1 \leq j \leq n-k} \left\{ \frac{\mathcal{W}_j^{(k)}}{\mathcal{W}_{j-1}^{(k)}} \right\} = \frac{\mathcal{W}_{n-k}^{(k)}}{\mathcal{W}_{n-k-1}^{(k)}}.$$

Lemma 3 is proved.

LEMMA 4. *Let* $\alpha = \min_{1 \leq j \leq n-k} \left\{ \dfrac{\mathcal{W}_j^{(k)}}{\mathcal{W}_{j-1}^{(k)}} \right\}$ *and*

$$\beta = \min \left\{ \mathcal{W}_1^{(k)},\ \min_{2 \leq j \leq n-k} \left\{ \frac{\mathcal{W}_j^{(k)} - \mathcal{W}_{j-1}^{(k)}}{\mathcal{W}_{j-1}^{(k)} - \mathcal{W}_{j-2}^{(k)}} \right\} \right\}.$$

*Then* $\alpha \geq \beta > \alpha - 1$.

*Proof.* Write $a_j = \dfrac{\mathcal{W}_j^{(k)}}{\mathcal{W}_{j-1}^{(k)}} \geq 0$ for $1 \leq j \leq n-k$ and $b_1 = \mathcal{W}_1^{(k)} \geq 0$, $b_j = \dfrac{\mathcal{W}_j^{(k)} - \mathcal{W}_{j-1}^{(k)}}{\mathcal{W}_{j-1}^{(k)} - \mathcal{W}_{j-2}^{(k)}} \geq 0$ for $2 \leq j \leq n-k$. Now, $a_1 = b_1$ and $a_j \geq b_j$ since by Lemma 3

$$\mathcal{W}_{j-1}^{(k)} \mathcal{W}_j^{(k)} - \mathcal{W}_j^{(k)} \mathcal{W}_{j-2}^{(k)} \geq \mathcal{W}_{j-1}^{(k)} \mathcal{W}_j^{(k)} - \left(\mathcal{W}_{j-1}^{(k)}\right)^2.$$

Therefore, $\min_{1 \leq j \leq n-k} \{a_j\} \geq \min_{1 \leq j \leq n-k} \{b_j\}$. Write

$$\alpha = \min_{1 \leq j \leq n-k} \{a_j\} \quad \text{and} \quad \beta = \min_{1 \leq j \leq n-k} \{b_j\}.$$

We know, by Lemma 3, that $\alpha = \dfrac{\mathcal{W}_{n-k}^{(k)}}{\mathcal{W}_{n-k-1}^{(k)}}$. Let us say

$$\beta = \frac{\mathcal{W}_m^{(k)} - \mathcal{W}_{m-1}^{(k)}}{\mathcal{W}_{m-1}^{(k)} - \mathcal{W}_{m-2}^{(k)}}$$

for some $2 \leq m \leq n-k$. Then,

$$1 \geq \frac{\beta}{\alpha} = \frac{\mathcal{W}_{n-k-1}^{(k)}}{\mathcal{W}_{n-k}^{(k)}} \cdot \frac{\mathcal{W}_m^{(k)} - \mathcal{W}_{m-1}^{(k)}}{\mathcal{W}_{m-1}^{(k)} - \mathcal{W}_{m-2}^{(k)}} = \frac{\mathcal{W}_{n-k-1}^{(k)}}{\mathcal{W}_{n-k}^{(k)}} \frac{\mathcal{W}_m^{(k)}}{\mathcal{W}_{m-1}^{(k)}} \left( \frac{1 - \dfrac{\mathcal{W}_{m-1}^{(k)}}{\mathcal{W}_m^{(k)}}}{1 - \dfrac{\mathcal{W}_{m-2}^{(k)}}{\mathcal{W}_{m-1}^{(k)}}} \right).$$

By $(*)$ of Lemma 3, $\dfrac{\mathcal{W}^{(k)}_{n-k-1}}{\mathcal{W}^{(k)}_{n-k}} \dfrac{\mathcal{W}^{(k)}_m}{\mathcal{W}^{(k)}_{m-1}} \geq 1$ and hence, by using $(*)$ again

$$\frac{\beta}{\alpha} \geq \frac{1 - \dfrac{\mathcal{W}^{(k)}_{m-1}}{\mathcal{W}^{(k)}_m}}{1 - \dfrac{\mathcal{W}^{(k)}_{m-2}}{\mathcal{W}^{(k)}_{m-1}}} \geq \frac{1 - \dfrac{\mathcal{W}^{(k)}_{n-k-1}}{\mathcal{W}^{(k)}_{n-k}}}{1 - \dfrac{1}{\mathcal{W}^{(k)}_1}} > 1 - \frac{\mathcal{W}^{(k)}_{n-k-1}}{\mathcal{W}^{(k)}_{n-k}} > 1 - \frac{1}{\alpha}.$$

Therefore $\alpha \geq \beta > \alpha - 1$ and the lemma is proved. □

Theorem II follows, since by Theorem I, $\mathcal{N}\bigl(i^*\mathcal{F}^d, \mathbf{V}^{n-k}_{(d_1,\ldots,d_k)}\bigr) > 0$, which happens, by Lemma 2, for $d \geq \beta$; and as, by Lemma 4, $\beta > \alpha - 1$ we obtain $d > \alpha - 1$. Now, $d$ is a positive integer and this gives $d \geq \alpha$. By the corollary to Theorem I, $\varrho_j\bigl(\mathbf{V}^{n-k}_{(d_1,\ldots,d_k)}\bigr) = (d_1 \cdots d_k)\mathcal{W}^{(k)}_j(d_1 - 1, \ldots, d_k - 1)$ so that

$$\alpha = \frac{\varrho_{n-k}\bigl(\mathbf{V}^{n-k}_{(d_1,\ldots,d_k)}\bigr)}{\varrho_{n-k-1}\bigl(\mathbf{V}^{n-k}_{(d_1,\ldots,d_k)}\bigr)}.$$

*Remark* 1. It is clear from the proof of Theorem II that all that is needed to obtain the bound $d \geq \dfrac{\varrho_{n-k}\bigl(\mathbf{V}^{n-k}_{(d_1,\ldots,d_k)}\bigr)}{\varrho_{n-k-1}\bigl(\mathbf{V}^{n-k}_{(d_1,\ldots,d_k)}\bigr)}$ are the following relations involving polar classes: $\varrho_1\varrho_{j-1} - \varrho_j \geq \varrho_1\varrho_{j-2} - \varrho_{j-1}$ and $\varrho_j^2 \geq \varrho_{j-1}\varrho_{j+1}$ for $1 \leq j \leq n-k$. It would be interesting to know if such relations hold for the polar classes of a variety $\mathbf{V}^{n-k} \xrightarrow{i} \mathbb{P}^n_{\mathbb{C}}$ which is not necessarily a complete intersection.

*Remark* 2. If $n-k$ is even, $\mathcal{N}\bigl(i^*\mathcal{F}^d, \mathbf{V}^{n-k}_{(d_1,\ldots,d_k)}\bigr)$ is automatically positive so, assuming the foliation is nondegenerate just along the variety, we cannot use the same arguments as given in Theorem II to relate $d$ to the polar classes of $\mathbf{V}^{n-k}_{(d_1,\ldots,d_k)}$. In [14] we considered the codimension 1 case, regardless of $n-1$ been even or odd, but assumed the foliation was nondegenerate in the whole of $\mathbb{P}^n_{\mathbb{C}}$. This allowed us to bound from above the number of singularities of the foliation in $\mathbf{V}^{n-1}_{d_1}$ by $d^n + d^{n-1} + \cdots + d + 1$, the total number of singularities of $\mathcal{F}^d$ in $\mathbb{P}^n_{\mathbb{C}}$, whenever $n-1$ is even. However, if $n-k$ is even and we make the stronger hypothesis that both $\mathbf{V}^{n-k}_{(d_1,\ldots,d_k)}$ and $\mathbf{V}^{n-k+1}_{(d_1,\ldots,d_{k-1})} \supset \mathbf{V}^{n-k}_{(d_1,\ldots,d_k)}$ are invariant by $\mathcal{F}^d$, which is also nondegenerate along $\mathbf{V}^{n-k+1}_{(d_1,\ldots,d_{k-1})}$ then, using the relation

$$\mathcal{W}^{(k)}_\delta(x_1, \ldots, x_k) = \mathcal{W}^{(k+1)}_\delta(x_1, \ldots, x_{k+1}) - x_{k+1}\mathcal{W}^{(k+1)}_{\delta-1}(x_1, \ldots, x_{k+1}),$$

the same argument given in the proof of Theorem II works, only this time we bound $\mathcal{N}(i^*\mathcal{F}^d, \mathbf{V}^{n-k}_{(d_1,\ldots,d_k)})$ from above by the corresponding number of singularities of $\mathcal{F}^d$ in $\mathbf{V}^{n-k+1}_{(d_1,\ldots,d_{k-1})}$. In this case we obtain precisely the same relations as in the odd dimensional case.

*Acknowlegements.* I am grateful to A. Lins Neto, I. Vainsencher and to R. S. Mol for pointing out mistakes in an earlier version of this work, to PRONEX - Dynamical Systems (Brasil), to Univ. de Rennes I (France) and to Univ. de Valladolid (Spain) for support and hospitality.

INSTITUTO DE CIÊNCIAS EXATAS, UNIVERSIDADE FEDERAL DE MINAS GERAIS, BELO HORIZONTE, BRAZIL
*E-mail address*: msoares@mat.ufmg.br